\documentclass[12pt,a4paper]{article}

\usepackage{latexsym,amsmath,amssymb,cite}

\newcommand{\proof}{\noindent{\bf Proof.}\ \ }
\newcommand{\qed}{\hfill $\blacksquare$}

\newcommand{\hh}{\rightarrow}
\newcommand{\hc}{\rightsquigarrow}

\newtheorem{theorem}{Theorem}[section]
\newtheorem{lemma}[theorem]{Lemma}
\newtheorem{conjecture}{Conjecture}

\begin{document}

\baselineskip=19pt
\parindent=0.5cm

\title{\LARGE{\bf An improved upper bound on the\\ 
adjacent vertex distinguishing\\
chromatic index of a graph}}

\author{Lianzhu Zhang
\thanks{Research supported  by NSFC (No. 11171279 
and No. 10831001)}\\
\normalsize School of Mathematical Science, Xiamen University, 
Xiamen 361005, China\\
\normalsize {\tt Email:zhanglz@xmu.edu.cn}\\
\and
Weifan Wang
\thanks{Research supported by NSFC (No. 11071223) 
and ZJNSFC (No. Z6090150)}\\
\normalsize Department of Mathematics, Zhejiang Normal University, 
Jinhua 321004, China\\
\normalsize {\tt Email:wwf@zjnu.cn}\\
\and
Ko-Wei Lih
\thanks{Research supported by NSC (No. 100-2517-S-001-001-MY3)}\\
\normalsize Institute of Mathematics, Academia Sinica, 
Taipei 10617, Taiwan\\
\normalsize   {\tt Email:makwlih@sinica.edu.tw}}

\date{}%\small August 9, 2012}
\maketitle

%%%%%%%%%%%%%%%%%%%%%%%%%%%%%%%%%%%%%%%%%%%%%%%%%%%%%%%%%%%%%%%%%%
%
\begin{abstract}
\noindent
An adjacent vertex distinguishing coloring of a graph $G$ is a proper 
edge coloring of $G$ such that any pair of adjacent vertices are 
incident with distinct sets of colors. The minimum number of colors 
needed for an adjacent vertex distinguishing coloring of $G$ is 
denoted by $\chi'_a(G)$. In this paper, we prove that $\chi_a'(G)
\leqslant \frac 52(\Delta+2)$ for any graph $G$ having maximum degree 
$\Delta$ and no isolated edges. This improves a  result in [S. Akbari, 
H. Bidkhori, N. Nosrati, $r$-Strong edge colorings of graphs,  Discrete
Math.  306 (2006), 3005-3010], which states that $\chi_a'(G)\leqslant
3\Delta$ for any graph $G$ without isolated edges.

\bigskip


\end{abstract}
%
%%%%%%%%%%%%%%%%%%%%%%%%%%%%%%%%%%%%%%%%%%%%%%%%%%%%%%%%%%%%%%%%%%

%%%%%%%%%%%%%%%%%%%%%%%%%%%%%%%%%%%%%%%%%%%%%%%%%%%%%%%%%%%%%%%%%%
%
\section{Introduction}
%
%%%%%%%%%%%%%%%%%%%%%%%%%%%%%%%%%%%%%%%%%%%%%%%%%%%%%%%%%%%%%%%%%%

All graphs considered in this paper are finite and without self-loops 
or multiple edges. In order to avoid trivialities, we also assume that
every graph has no isolated vertices. Let $V(G)$ and $E(G)$ denote the 
vertex and the edge sets of $G$, respectively. Let $N_G(v)$ denote the 
set of neighbors of $v$ in $G$ and $d_G(v)=|N_G(v)|$ the degree of $v$ 
in $G$. A vertex $v$ is called a {\em $k$-vertex} if $d_G(v)=k$. Let 
$\Delta(G)$ and $\delta(G)$ denote the maximum and minimum degree of a 
vertex in $G$, respectively. An {\em edge $k$-coloring} of a graph $G$ 
is a function $\phi: E(G) \to \{1,2,\ldots ,k\}$ such that any two
incident edges receive different colors. The {\em  chromatic index},
denoted by $\chi'(G)$, of a graph $G$ is the smallest integer $k$ such
that $G$ has an edge $k$-coloring. Given an edge $k$-coloring $\phi$ of 
$G$, we use $C_{\phi}(v)$ to denote the set of colors assigned to edges 
incident to a vertex $v$. We call $C_{\phi}=\cup_{v \in V(G)} C_{\phi}(v)$ 
the color set of $\phi$. The coloring $\phi$ is called an {\em adjacent 
vertex distinguishing edge coloring} if $C_{\phi}(u)\neq C_{\phi}(v)$ for 
any pair of adjacent vertices $u$ and $v$. A graph $G$ is {\em normal} if 
it contains no isolated edges. Clearly, $G$ has an adjacent vertex 
distinguishing edge coloring if and only if $G$ is normal. The {\em 
adjacent vertex distinguishing chromatic index} $\chi'_{a}(G)$ of a graph 
$G$ is the smallest integer $k$ such that $G$ has an adjacent vertex 
distinguishing edge $k$-coloring.

Zhang, Liu and Wang \cite{Zhang02} first introduced and investigated
the adjacent vertex distinguishing edge coloring ({\em adjacent strong 
edge coloring} in their terminology) of graphs. They proposed the 
following conjecture.

\begin{conjecture}\label{conj}  If a connected normal graph $G$ is
different from a 5-cycle and  satisfies $|V(G)|\geqslant 3$, then 
$\chi'_{a}(G)\leqslant\Delta(G)+2$.
\end{conjecture}

Balister et al. \cite{Balister07} confirmed Conjecture \ref{conj} for 
all normal graphs $G$ that are bipartite or satisfy $\Delta(G)=3$. 
In particular, we need the following statement in the sequel.

\begin{theorem}\label{deta3}
For any normal graph $G$ with $\Delta(G)\leqslant 3$, $ \chi_a'(G)
\leqslant 5$.
\end{theorem}

They further proved that $\chi'_{a}(G)\leqslant \Delta(G)+O(\log k)$, 
where $k$ is the (vertex) chromatic number of the normal graph $G$. 
It follows from Brooks' Theorem that $\chi'_{a}(G)\leqslant 2\Delta(G)$
for $G$ with sufficiently large $\Delta(G)$. Hatami \cite{Hatami05} 
showed that every normal graph $G$ with $\Delta(G)>10^{20}$ has 
$\chi'_{a}(G)\leqslant \Delta(G)+300$ by the probabilistic method. 
Edwards et al. \cite{edwa} proved that $\chi'_a(G)\leqslant\Delta(G)+1$ 
if $G$ is a planar bipartite normal graph with $\Delta(G)\geqslant 12$. 
Wang and Wang \cite{Wang10} verified Conjecture \ref{conj} for a class 
of graphs with small maximum average degree. Their results were further 
extended by Hocquard and Montassier \cite{Hocquard11, Hocquard12}. 
Recently, it has been characterized in \cite{Wang11} which of the two 
cases $\chi'_{a}(G)=\Delta(G)$ and $\chi'_{a}(G)=\Delta(G)+1$ holds for 
a $K_4$-minor-free normal graph $G$ with $\Delta(G) \geqslant5$.

An adjacent vertex distinguishing edge coloring of a graph $G$ is a 
special case of a {\em vertex distinguishing edge coloring}, which 
requires that every pair of vertices be incident with distinct color 
sets. This more general notion was introduced by Burris and Schelp 
\cite{BS97}, and independently by Hor\v{n}\'{a}k and Sot\'{a}k \cite{HS95}, 
and \v{C}ern\'{y} et al. \cite{CHS96} (under the name {\em observability}). 
The reader is referred to \cite{Balister01, BBS02, BKLS04, BRS03, BHLW99, 
BHLW01, LL10a} for relevant results.

The aim of this paper is to improve the following upper bound obtained in 
\cite{akb}.

\begin{theorem}\label{deta5}
For any normal graph $G$, $\chi_a'(G)\leqslant 3\Delta(G)$.
\end{theorem}

The proof of our main theorem in Section 2 is based on an edge-partition 
result. The details will be supplied in the last section. In Section 3, 
the new upper bound is further reduced for regular graphs.

%%%%%%%%%%%%%%%%%%%%%%%%%%%%%%%%%%%%%%%%%%%%%%%%%%%%%%%%%%%%%%%%%%
%
\section{An improved upper bound}
%
%%%%%%%%%%%%%%%%%%%%%%%%%%%%%%%%%%%%%%%%%%%%%%%%%%%%%%%%%%%%%%%%%%

For a graph $G$ and any $S \subseteq E(G)$, the {\em edge-induced} 
subgraph $G[S]$ is the subgraph of $G$ whose edge set is $S$ and 
whose vertex set consists of all end vertices of edges in $S$. We 
only deal with subgraphs that are edge-induced subgraphs unless 
otherwise stated. For a subgraph $H$ of $G$, we use $\overline{H}$ 
to denote the edge-induced subgraph $G[E(G) \setminus E(H)]$ and 
call it the {\em complement} of $H$ in $G$. An {\em edge-partition} 
of a graph $G$ into subgraphs $G_1,G_2,\ldots ,G_m$ is a decomposition 
of $G$ that satisfies $V(G)=\cup_{i=1}^{m} V(G_i)$, $E(G)=\cup_{i=1}^{m} 
E(G_i)$ and $E(G_i)\cap E(G_j)=\emptyset$ for any pair $i\ne j$. 
Clearly, a subgraph $H$ of $G$ together with its complement $\overline{H}$ 
constitute an edge-partition of $G$. This edge-partition is said to be 
induced by the subgraph $H$. The proof of the following is deferred to 
Section 4.

\begin{theorem}\label{p1}
Let $G$ be a normal graph with $\Delta(G)\geqslant 6$. Then there is an 
edge-partition of $G$ induced by a subgraph $H$ such that the following 
conditions hold.
\begin{enumerate}
\item
Both $H$ and $\overline{H}$ are normal.
\item
$\Delta(H)\leqslant 3$.
\item
$\Delta(\overline{H}) \leqslant \Delta(G)-2$.
\end{enumerate}
\end{theorem}

\begin{theorem}\label{p2}
Let $G$ be a normal graph with $\Delta(G)\geqslant 4$. Then there is an 
edge-partition of $G$ into subgraphs $G_0, G_1, \ldots , G_k$, $k 
\leqslant \lfloor \Delta(G)/2 \rfloor-2$, such that the following hold.
\begin{enumerate}
\item
Every $G_i$ is a normal subgraph.
\item
$\Delta(G_i)\leqslant 3$ for $1 \leqslant i \leqslant k$.
\item
$\Delta(G_0)\leqslant 5$.
\end{enumerate}
\end{theorem}
 
\proof  
The proof proceeds by induction on $\Delta(G)$. If $\Delta(G)\leqslant 5$, 
the result holds trivially. Let $G$ be a normal graph with $\Delta(G)
\geqslant 6$. By Theorem \ref{p1}, there is an edge-partition of $G$ 
induced by a subgraph $H$ such that both $H$ and $\overline{H}$ are normal, 
$\Delta(H)\leqslant 3$ and $\Delta(\overline{H}) \leqslant \Delta(G)-2$. 
Clearly, $\Delta(\overline{H}) \geqslant 3$. If $\Delta (\overline{H})=3$, 
then $\Delta(G)=6$. Let $G_0=H$ and $G_1=\overline{H}$. If $\Delta(\overline{H})
\geqslant 4$, by the induction hypothesis, there is an edge-partition of 
$\overline{H}$ into subgraphs $G_0, G_1, \ldots , G_k$, $k\leqslant \lfloor 
\Delta(\overline{H})/2 \rfloor-2$, such that properties 1, 2 and 3 hold. Now 
let $G_{k+1}=H$. Then $G_0, G_1, \ldots , G_k, G_{k+1}$ form an edge-partition 
of $G$. Note that $k+1\leqslant  \lfloor \Delta(\overline{H})/2 \rfloor-2+1 
\leqslant \lfloor (\Delta(G)-2)/2 \rfloor-1 = \lfloor \Delta(G)/2 \rfloor-2$ 
and we are done. 
\qed

\begin{lemma}\label{sum}
If a normal graph $G$ has an edge-partition into two normal subgraphs
$G_1$ and $G_2$, then $\chi_a'(G)\leqslant \chi'_a(G_1)+\chi_a'(G_2).$
\end{lemma}

\proof
For $i=1,2$, let $\phi_i$ be an adjacent vertex distinguishing edge 
coloring of $G_i$ satisfying $|C_{\phi_i}|=\chi'_a(G_i)$ and $C_{\phi_1}
\cap C_{\phi_2}=\emptyset$. The union of $\phi_1$ and $\phi_2$ forms a 
proper edge coloring $\phi$ of $G$ with color set $C_{\phi_1}\cup C_{\phi_2}$. 
Let $uv\in E(G)$ with $d_G(u)=d_G(v)$. Since $E(G_1)\cap E(G_2)=\emptyset$, 
we may assume that $uv\in E(G_1)\setminus E(G_2)$ with $d_{G_1}(u) \geqslant 
d_{G_1}(v)$. Since $G_1$ is normal, $uv$ is not an isolated edge of $G_1$, 
i.e., $d_{G_1}(u)\geqslant 2$.  By definition of $\phi_1$, there exists a 
$c\in C_{\phi_1}(u)\setminus C_{\phi_1}(v)$.  Since $C_{\phi_1}\cap C_{\phi_2}
=\emptyset$, it follows that $c \in C_{\phi}(u) \setminus C_{\phi}(v)$, and 
hence $C_{\phi}(u) \ne C_{\phi}(v)$. Consequently, $\chi'_a(G) \leqslant 
|C_{\phi_1}\cup C_{\phi_2}|=|C_{\phi_1}|+|C_{\phi_2}|=\chi'_a(G_1)+\chi'_a(G_2)$.  
\qed

\begin{theorem}\label{main}
If $G$ is a normal graph, then $\chi'_a(G)\leqslant \frac 52(\Delta(G)+2)$.
\end{theorem}

\proof 
The result can be derived immediately from Theorem \ref{deta3} when 
$\Delta(G) \leqslant 3$. Now assume that $\Delta(G) \geqslant 4$. By Theorem 
\ref{p2}, there is an edge-partition of $G$ into subgraphs $G_0, G_1, 
\ldots , G_k$, $k\leqslant \lfloor \Delta(G)/2 \rfloor-2$, such that 
properties 1, 2 and 3 hold. Using Lemma \ref{sum} and Theorem \ref{deta3} 
repeatedly, we have
\begin{eqnarray*}
\chi_a'(G) & \leqslant & \chi_a'(G_0)+\chi_a'(G_1)+\cdots +\chi_a'(G_k)\\
           & \leqslant & \chi_a'(G_0)+5k\\
           & \leqslant & \chi_a'(G_0)+5(\lfloor \Delta(G)/2  \rfloor-2).
\end{eqnarray*}
By Theorem \ref{p2}, $\Delta(G_0)\leqslant 5$. It follows from Theorem 
\ref{deta5} that $\chi_a'(G)\leqslant 15+5(\lfloor \Delta(G)/2 \rfloor-2) 
\leqslant \frac 52(\Delta(G)+2)$. 
\qed

%%%%%%%%%%%%%%%%%%%%%%%%%%%%%%%%%%%%%%%%%%%%%%%%%%%%%%%%%%%%%%%%%%
%
\section{Regular graphs}
%
%%%%%%%%%%%%%%%%%%%%%%%%%%%%%%%%%%%%%%%%%%%%%%%%%%%%%%%%%%%%%%%%%%

Theorem \ref{main} can be further improved for regular graphs. 
We first establish an auxiliary edge-partition lemma. We need 
the following well-known result of Vizing \cite{vizi} on chromatic 
index.

\begin{theorem}\label{vizing}
For every graph $G$, $\chi'(G) \leqslant \Delta(G)+1$.
\end{theorem}

\begin{lemma}\label{regular}
Let $G$ be a regular graph of degree $r\geqslant 5$. Then there is an 
edge-partition of $G$ into normal subgraphs $G_1, G_2, \ldots , G_k$ 
such that one of the following conditions holds.
\begin{enumerate}
\item
If $r \equiv 2 \pmod{3}$, then $k=(r+1)/3$ and $\Delta(G_i) \leqslant 
3$ for $1 \leqslant i \leqslant k$.
\item
If $r \equiv 1 \pmod{3}$, then $k=(r-1)/3$, $\Delta(G_i) \leqslant 4$ 
for $1 \leqslant i \leqslant 2$ and $\Delta(G_i)\leqslant 3$ for $3 
\leqslant i \leqslant k$.
\item
If $r \equiv 0 \pmod{3}$, then $k=r/3$ and $\Delta(G_1) \leqslant 4$  
and $\Delta(G_i)\leqslant 3$ for $2 \leqslant i \leqslant k$.
\end{enumerate}
\end{lemma}

\proof  
By Theorem \ref{vizing}, $E(G)$ can be partitioned into $r+1$ disjoint 
color classes $E_1, E_2,\ldots , E_{r+1}$ such that each $E_i$ is a 
matching of $G$. Let $H$ be a subgraph of $G$ edge-induced by $m$, $3 
\leqslant m \leqslant r$, of these color classes. Obviously, $\Delta(H) 
\leqslant m$. For any given vertex $v$ of $G$, exactly one color is not 
used on any edge incident with $v$ since $G$ is $r$-regular. Therefore 
$d_H(v) \geqslant 2$, and hence $H$ is a normal graph.

If $r \equiv 2 \pmod{3}$, let $k= (r+1)/3$. Then we define $G_1=G[E_1 
\cup E_2 \cup E_3]$, $G_2=G[E_4 \cup E_5 \cup E_6]$, $\ldots ,$ $G_k=
G[E_{r-1} \cup E_{r} \cup E_{r+1}]$. Then $G_1, G_2, \ldots , G_k$ 
form an edge-partition of $G$ satisfying condition 1.

If $r \equiv 1 \pmod{3}$, let $k=(r-1)/3$. Then we define  $G_1=G[E_1 
\cup E_2 \cup E_3 \cup E_4]$, $G_2=G[E_5 \cup E_6 \cup E_7 \cup E_8]$, 
$G_3=[E_9 \cup E_{10} \cup E_{11}]$, $\ldots ,$ $G_k=G[E_{r-1} \cup E_{r} 
\cup E_{r+1}]$. Then $G_1, G_2, \ldots , G_k$ form an edge-partition of 
$G$ satisfying condition 2.

If $r \equiv 0 \pmod{3}$, let $k=r/3$. Then we define $G_1=G[E_1 \cup E_2 
\cup E_3 \cup E_4]$, $G_2=G[E_5 \cup E_6 \cup E_7]$, $G_3=[E_8 \cup E_9 
\cup E_{10}]$, $\ldots ,$ $G_k=G[E_{r-1} \cup E_{r} \cup E_{r+1}]$. Then 
$G_1, G_2, \ldots , G_k$ form an edge-partition of $G$ satisfying condition 3.
\qed

\begin{theorem}\label{5/3}
Let $G$ be a regular graph of degree $r\geqslant 2$. Then $\chi_a'(G)
\leqslant (5r + 37)/3.$
\end{theorem}

\proof
If $2 \leqslant r \leqslant 4$, the result follows from Theorems \ref{deta3} 
and \ref{deta5}. Assume that $r\geqslant 5$. By Lemma \ref{regular}, there 
is an edge-partition of $G$ into normal subgraphs $G_1, G_2, \ldots , G_k$
such that one of the stated conditions 1, 2 or 3 holds.

If condition 1 holds, by Lemma \ref{sum}, Theorems \ref{deta3} and
\ref{deta5}, we have $\chi'_a(G) \leqslant \sum_{i=1}^{k}\chi'_a(G_i)$ 
$\leqslant 5k= 5(r+1)/3 < (5r + 37)/3.$

If condition 2 holds, then $\chi'_a(G) \leqslant \chi_a'(G_1)+\chi'_a(G_2)+
\sum_{i=3}^{k}\chi'_a(G_i) \leqslant 12+12+ 5(k-2)= 5(r-1)/3 + 14 = 
(5r + 37)/3$.

If condition 3 holds, then $\chi'_a(G) \leqslant \chi_a'(G_1)+
\sum_{i=2}^{k}\chi'_a(G_i) \leqslant 12+5(k-1) = 5r/3 + 7 < (5r + 37)/3$.
\qed

\bigskip

Note that the upper bound in Theorem \ref{5/3} is better than the upper 
bound in Theorem \ref{main} when $r \geqslant 14$.

%%%%%%%%%%%%%%%%%%%%%%%%%%%%%%%%%%%%%%%%%%%%%%%%%%%%%%%%%%%%%%%%%%
%
\section{Proof of Theorem \ref{p1}}
%
%%%%%%%%%%%%%%%%%%%%%%%%%%%%%%%%%%%%%%%%%%%%%%%%%%%%%%%%%%%%%%%%%%

We devote this section to a complete proof of Theorem \ref{p1}.

Assume that $G$ is a normal graph with $\Delta(G) \geqslant 6$.
We abbreviate $\Delta(G)$ and $d_G(v)$ to $\Delta$ and $d(v)$, 
respectively. Let $\mathcal{H}(G)$ be the collection of subgraphs 
$M$ of $G$ that satisfy the following conditions.

1.\  $\Delta(M)\leqslant 3$.

2.\  If $d(v)=\Delta$, then $d_M(v)\geqslant 2$.

3.\  If $d(v)=\Delta - 1$, then $d_M(v)\geqslant 1$.

We first show that $\mathcal{H}(G) \ne \emptyset$. By Theorem 
\ref{vizing}, $E(G)$ can be partitioned into $\Delta + 1$ 
disjoint color classes $E_1, E_2,\ldots , E_{\Delta +1}$ such 
that each $E_i$ is a matching of $G$. Let $M=G[E_1\cup E_2\cup
E_3]$. Then $\Delta(M)\leqslant 3$.  For a $\Delta$-vertex $x$ 
of $G$, at most one among $E_1,E_2,E_3$ contains no edge incident 
with $x$. For a $(\Delta -1)$-vertex $y$ of $G$, at most two among 
$E_1,E_2,E_3$ contain no edge incident with $y$. Thus $M\in 
\mathcal{H}(G)$. 

For any $M\in \mathcal{H}(G)$, it is easy to see that $\Delta(\overline{M})
\leqslant \Delta-2$. Now let $I(M)$ and $I(\overline{M})$ denote
the sets of isolated edges of $M$ and $\overline{M}$, respectively, 
and write $i(M)=|I(M)|$ and $i(\overline{M})=|I(\overline{M})|$.
Among all subgraphs $M$ that attain the minimum for $i(M)+i(\overline{M})$,
we pick and fix an $H$ that has minimum number of edges.

We are going to show that the edge-partition of $G$ induced by this $H$ 
satisfies conditions 1, 2 and 3 of Theorem \ref{p1}. If $i(H)+
i(\overline{H})=0$, then we are done. Now we assume that $i(H)+
i(\overline{H})>0$. 

\bigskip

We first classify some of the vertices of $G$ into two types.

A vertex $v\in V(G)$ is classified as {\em type-I} if $1\leqslant 
d_H(v)\leqslant 2$, $d(v)\geqslant \Delta-1$, and for every $u\in 
N_{\overline{H}}(v)$, one of the following three conditions holds.

{\rm (1)}\  $d_H(u)=3$.

{\rm (2)}\  $d_H(u)=d_{\overline{H}}(u)=2$.

{\rm (3)}\  $d_H(u)\leqslant 1$, $d_{\overline{H}}(u)=2$, and, for 
the unique $w\in N_{\overline{H}}(u)\setminus \{v\}$, both 
$d_{\overline{H}}(w)=1$ and $d_H(w)=3$.

\bigskip

\noindent 
{\bf Claim 1.}\ 
{\em Suppose that $vv'\in I(H)$ with $d(v)\geqslant d(v')$. Then 
$d(v)=\Delta-1$ and $v$ is a type-I vertex.}

\bigskip

\proof 
Since $H\in \mathcal{H}(G)$ and $vv'$ is an isolated edge of $H$, 
$d_H(v)=1$ and $d(v)\leqslant \Delta-1$. If $d(v)\leqslant \Delta-2$, 
then $H'=H \setminus \{vv'\}\in \mathcal{H}(G)$. Note that $i(H')=
i(H)-1$ and $i(\overline{H'})\leqslant i(\overline{H})$ since $vv' 
\not \in I(\overline{H'})$. The subgraph $H'$ contradicts the choice 
of $H$. Consequently, $d(v)=\Delta-1$.

Assume to the contrary that $v$ is not a type-I vertex. Then there 
exists a particular $u \in N_{\overline{H}}(v)$ that satisfies none 
of (1), (2) or (3). Thus, the following three statements hold for this $u$.

(a)\  $d_H(u) \ne 3$, and hence $d_H(u) \leqslant 2$.

(b)\  If $d_H(u)=2$, then $d_{\overline{H}}(u) \ne 2$.

(c)\  If $d_H(u)\leqslant 1$ and $d_{\overline{H}}(u)=2$, then, for the 
unique $w\in N_{\overline{H}}(u)\setminus \{v\}$, $d_{\overline{H}}(w)
=1$ implies $d_H(w) \ne 3$, and hence $d_H(w) \leqslant 2$.

Define $H'=H\cup \{uv\}$ for case (b) or when $d_{\overline{H}}(w) 
\ne 1$ for case (c). Define $H'=H\cup \{uv,uw\}$ when $d_{\overline{H}}(w)
= 1$ for case (c). It is easy to check that $H'\in \mathcal{H}(G)$. Since 
$d_{\overline{H'}}(v) = d(v) - d_{H'}(v) = (\Delta - 1)-2 > 2$, no new 
isolated edge is created in $\overline{H'}$. Yet $i(H')=i(H)-1$. This 
contradicts the choice of $H$.
\qed

\bigskip

A vertex $u\in V(G)$ is classified as {\em type-II} if $d_H(u)=3$, or 
$d_H(u)=d_{\overline{H}}(u)=2$, and for every $v\in N_H(u)$, one of the 
following two conditions holds.

{\rm (4)}\ $1\leqslant d_H(v)\leqslant 2$ and $d(v)\geqslant \Delta-1$.

{\rm (5)}\  $d_H(v)= 2$, $d(v)<\Delta-1$, and, for the unique $w\in N_H(v)
\setminus \{u\}$, both $d_H(w)=1$ and $d(w)=\Delta-1$.

\bigskip

\noindent 
{\bf Claim 2.}\ 
{\em Suppose that $uu'\in I(\overline{H})$ with $d(u)\geqslant d(u')$. 
Then $d_H(u)=3$ and $u$ is a type-II vertex.}

\bigskip

\proof  
Since $uu'$ is an isolated edge of $\overline{H}$ and $G$ has no 
isolated edges, it follows that $d_H(u)\geqslant 1$. If $d_H(u)
\leqslant 2$, then  $H'=H \cup \{uu'\} \in \mathcal{H}(G)$. Note 
that $i(H') \leqslant i(H)$ and $i(\overline{H'})= i(\overline{H})-1$. 
The subgraph $H'$ contradicts the choice of $H$. Consequently, $d_H(u)=3$.

Assume to the contrary that $u$ is not a type-II vertex. Then there 
exists a particular $v \in N_H(u)$ that satisfies neither (4) nor (5). 
Thus, the following two statements hold for this $v$.

(d)\  If $1\leqslant d_H(v)\leqslant 2$, then $d(v)< \Delta -1$.

(e)\  If $d_H(v)= 2$, $d(v)<\Delta-1$, then, for the unique $w \in 
N_H(v)\setminus \{u\}$, $d_{H}(w)=1$ implies $d(w) \ne \Delta-1$, and 
hence $d(w) < \Delta -1$.

If $d_H(v)=1$ or $d_H(v)=2$ and $d_H(w) \geqslant 2$, let $H'=H 
\setminus \{uv\}$. If $d_{H}(v)= 2$ and $d_H(w)=1$, let $H'=H 
\setminus \{uv,vw\}$. Thus, the subgraph $H' \in \mathcal{H}(G)$ 
and satisfies $i(H')\leqslant i(H)$ and $i(\overline{H'})=
i(\overline{H})-1$, contradicting the choice of $H$. 
\qed

\bigskip

We observe that no vertex can be classified both as type-I and 
type-II since $1\leqslant d_H(z)\leqslant 2$ and $d(z)\geqslant 
\Delta-1\geqslant 5$ for a type-I vertex $z$, while $d_H(w)= 3$ 
or $d_H(w)=d_{\overline{H}}(w)=2$ for a type-II vertex $w$. 

An {\em $H$-chain} emanating from a vertex $u$ is a path from $u$ 
to a $v\in N_H(u)$ when $v$ satisfies (4), or through $v$ to the 
unique $w\in N_H(v)\setminus \{u\}$ when $v$ satisfies (5). We write 
$u \hh x$ for an $H$-chain emanating from $u$ and terminating at $x$.
An {\em $\overline{H}$-chain} emanating from a vertex $v$ is a path 
from $v$ to a $u\in N_{\overline{H}}(v)$ when $u$ satisfies (1) or (2), 
or through $u$ to the unique $w\in N_{\overline{H}}(u) \setminus \{v\}$ 
when $u$ satisfies (3). We write $v \hc y$ for an $\overline{H}$-chain 
emanating from $v$ and terminating at $y$. A path $P$ of $G$ is called 
an {\em alternating chain} if $P$ is a concatenation of $H$-chains and 
$\overline{H}$-chains such that they appear alternately and the 
terminating vertex of one chain is the emanating vertex of the next chain.

\bigskip

\noindent 
{\bf Claim 3.}\ 
{\em  If $vv'\in I(H)$ satisfies $d(v)\geqslant d(v')$, then
the two ends of each $H$-chain or $\overline{H}$-chain of an
alternating chain $P$ beginning with $v$ are of different types.}

\bigskip

\proof  
Let $v_0=v$. By Claim 1, $v_0$ is a type-I vertex. By the definition 
of an alternating chain, we may assume that $P$ is $v_0 \hc u_1 \hh 
v_1 \hc \cdots \hh v_{s-1} \hc u_s$ or $P$ is $v_0 \hc u_1 \hh v_1 
\hc \cdots  \hc u_s \hh v_s$, where $s\geqslant 1$. It suffices to 
prove by induction that $v_1, v_2, \ldots , v_s$ are type-I vertices 
and $u_1, u_2, \ldots , u_s$ are type-II vertices. Equivalently, for 
each $1 \leqslant k\leqslant s$, the following statements (A) and (B) 
are true.

(A)\  If $v_1, v_2, \ldots , v_{k-1}$ are type-I vertices and $u_1, 
u_2, \ldots , u_{k-1}$ are type-II vertices, then $u_k$ is a type-II 
vertex.

(B)\  If $v_1, v_2, \ldots , v_{k-1}$ are type-I vertices and $u_1, 
u_2, \ldots , u_k$ are type-II vertices, then $v_k$ is a type-I vertex.

In order to show (A), assume to the contrary that $u_k$ is not a
type-II vertex.  Since $v_{k-1} \hc u_k$ and $v_{k-1}$ is a type-I
vertex, $d_H(u_k)=3$, or $d_H(u_k)=d_{\overline{H}}(u_k)=2$. Then 
there exists a vertex $x\in N_H(u_k)$ such that the following two 
statements hold for this $x$.

(d$'$)\  If $1\leqslant d_H(x)\leqslant 2$, then $d(x)< \Delta -1$.

(e$'$)\  If $d_H(x)= 2$, $d(x)<\Delta-1$, then, for the unique $y \in 
N_H(x)\setminus \{u_k\}$, $d_H(y)=1$ implies $d(y) < \Delta-1$.

Since $v_0, v_1, \ldots , v_{k-1}$ are type-I vertices by the induction 
hypothesis, $1\leqslant d_H(v_i)\leqslant 2$ and $d(v_i)\geqslant 
\Delta-1$ for all $0 \leqslant i \leqslant k-1$. Since $d_H(x)=3$, or 
$d(x)<\Delta-1$, it follows that $x \notin \{v_0, v_1, \ldots , 
v_{k-1}\}$. We next show that $x\notin \{u_1,u_2,\ldots,u_{k-1}\}$.

Assume to the contrary that there is an index $i$ ($i<k$) such that 
$x=u_i$. Since $u_i$ is a type-II vertex and $u_k\in N_H(u_i)$, it 
follows that $d_H(u_k)\leqslant 2$. We have already known that 
$d_H(u_k)=3$, or $d_H(u_k)=d_{\overline{H}}(u_k)=2$. Hence, $d_H(u_k)=2$ 
and $d(u_k)=4$. Let $z\in N_H(u_k)\setminus \{u_i\}$. Define
\[
H'=(H \cup \bigcup_{j=0}^{i-1}E(v_j \hc u_{j+1}))\setminus
(S\cup \bigcup_{j=1}^{i-1}E(u_j \hh v_{j})), 
\]
where $S=\{u_iu_k,u_kz\}$ if $d_H(z)=1$; or $S=\{u_iu_k\}$ otherwise.  
It is straightforward to check that $H'\in \mathcal{H}(G)$ such that
$i(H')=i(H)-1$ and $i(\overline{H'}) = i(\overline{H})$, which
contradicts the choice of $H$. 

Suppose that $d_H(x)=1$ or $d_H(x)= 2$ and $d_H(y)>1$ in (e$'$). If 
$d_H(u_k)=3$, then let $H'=H \setminus \{xu_k\}$. It is obvious that 
$H' \in \mathcal{H}(G)$. Since $xu_k$ is adjacent to an edge in $v_{k-1} 
\hc u_k$, $xu_k$ can not be an isolated edge of $\overline{H'}$. Thus, 
$i(H')=i(H)$ and $i(\overline{H'})=i(\overline{H})$. However, $|E(H')|=
|E(H)|-1$, which contradicts the choice of $H$. If $d_H(u_k)=
d_{\overline{H}}(u_k)=2$, define
\[
H'=(H\cup
 \bigcup_{i=0}^{k-1}E(v_i\hc u_{i+1}))\setminus
 (\bigcup_{i=1}^{k-1}E(u_i\hh v_i)  \cup \{ xu_k\}).
\]
Note that $d_{H'}(u_i)=d_H(u_i)$ and $d_{H'}(v_i)=d_H(v_i)$ for $1 
\leqslant i \leqslant k$, $d_{H'}(v_0)=d_H(v_0)+1=2$, $d_{\overline{H'}}(v_0)
=(\Delta-1)-2\geqslant 3$, and hence $v'v_0\notin I(H')$. It follows that 
$i(H')=i(H)-1$ and $i(\overline{H'})=i(\overline{H})$, which contradicts 
the choice of $H$.

Next consider the case $d_H(y)=1$ in (e$'$). Then $y \notin \{v_0, v_1, 
\ldots , v_{k-1}\}$ since $d(y) < \Delta-1$; $y\notin \{u_1,u_2,\ldots,
u_{k-1}\}$ for each type-II vertex $u_i$ ($1 \leqslant i \leqslant k-1$) 
has $d_H(u_i)\geqslant 2$. Define
\[
H'=(H\cup
 \bigcup_{i=0}^{k-1}E(v_i \hc u_{i+1}))\setminus
 (\bigcup_{i=1}^{k-1}E(u_i \hh v_i)  \cup \{ xy, xu_k\}).
\]
Then $H'\in \mathcal{H}(G)$. Reasoning as before, we see that $i(H')=
i(H)-1$ and $i(\overline{H'})=i(\overline{H})$, which contradicts the 
choice of $H$.

To prove (B), assume to the contrary that $v_k$ is not a type-I vertex.  
Since $u_k \hh v_k$ and $u_k$ is a type-II vertex, $1 \leqslant d_H(v_k) 
\leqslant 2$ and $d(v_k) \geqslant \Delta -1$. Then there exists a vertex
$x\in N_{\overline{H}}(v_k)$ such that the following three statements hold 
for this $x$.

(a$'$)\  $d_H(x) \ne 3$, and hence $d_H(x) \leqslant 2$.

(b$'$)\  If $d_H(x)=2$, then $d_{\overline{H}}(x) \ne 2$.

(c$'$)\  If $d_H(x)\leqslant 1$ and $d_{\overline{H}}(x)=2$,  then, for  
the unique $y\in N_{\overline{H}}(x)\setminus \{v_k\}$, $d_{\overline{H}}(y)
=1$ implies $d_H(y) \leqslant 2$.

Since $u_1, u_2, \ldots , u_k$ are type-II vertices by the induction
hypothesis,  we see that for $1\leqslant i\leqslant k$,  either 
$d_H(u_i)=3$ or $d_H(u_i)=d_{\overline{H}}(u_i)=2$. Therefore, $x \notin 
\{u_1, u_2, \ldots , u_k\}$.

We next show that $x\notin \{v_0, v_1, \ldots , v_{k-1} \}$. Assume to 
the contrary that there is an index $i$ ($0 \leqslant i \leqslant k-1$) 
such that $x=v_i$. Since $v_i$ is a type-I vertex and $v_k \in 
N_{\overline{H}}(v_i)$, it follows that $d_H(v_k)=3$ or $d_H(v_k)=
d_{\overline{H}}(v_k)=2$. However, $d_H(v_k)\leqslant 2$ and $d(v_k)
\geqslant \Delta-1\geqslant 5$ since $u_k \hh v_k$. We have reached a 
contradiction.

Now assume $d_{\overline{H}}(y)=1$ in (c$'$). Then $y \notin \{u_1, 
u_2, \ldots , u_k\}$. We also have $y\notin \{v,v_1,\ldots,u_{k-1}\}$, 
for otherwise it would imply $d_{\overline{H}}(y)\geqslant 2$. Define
\[
H'=(H\cup S \cup
 \bigcup_{i=0}^{k-1}E(v_i\hc u_{i+1}))\setminus
 \bigcup_{i=1}^{k-1}E(u_i \hh v_i),
\]
where $S=\{xy,xv_k\}$ when $d_{\overline{H}}(y)=1$ for case (c$'$); 
$S=\{xv_k\}$ for case (b$'$) or when $d_{\overline{H}}(y)\ne 1$ for 
case (c$'$). It is easy to check that $H'\in \mathcal{H}(G)$ such 
that $i(H')=i(H)-1$ and $i(\overline{H'})=i(\overline{H})$. This 
contradicts the choice of $H$. 
\qed

\bigskip

\noindent 
{\bf Claim 4.}\ 
{\em  If $uu'\in I(\overline{H})$ satisfies $d(u)\geqslant d(u')$, then
the two ends of each $H$-chain or $\overline{H}$-chain of an alternating
chain $P$ beginning with $u$ are of different types.}

\bigskip

\proof  
Let $u_1=u$ which is a type-II vertex by Claim 2. By the definition 
of an alternating chain, we may assume that $P$ is $u_1\hh v_1\hc u_2 
\hh \cdots \hc u_s \hh v_s$ or $P$ is $u_1\hh v_1 \hc u_2 \hh \cdots 
\hh v_{s-1}\hc u_s$,  where $s\geqslant 1$. Similar to the proof of 
Claim 3, we may argue that, for each $1\leqslant k\leqslant s$, the 
following statements (C) and (D) are true.

(C)\  If $u_1, u_2 \ldots , u_k$ are type-II vertices and $v_1, v_2, 
\ldots , v_{k-1}$  are type-I vertices, then $v_k$ is a type-I vertex.

(D)\ If $u_1, u_2 \ldots , u_{k-1}$ are type-II vertices and $v_1, v_2, 
\ldots , v_{k-1}$  are type-I vertices, then $u_k$ is a type-II vertex.

The proof of (B) in Claim 3 can be adapted to show the validity of (C). 
Here we define
\[
H'=(H\cup S \cup
 \bigcup_{i=1}^{k-1}E(v_i\hc u_{i+1}))\setminus
 \bigcup_{i=1}^{k-1}E(u_i \hh v_i),
\]
where $S=\{xy,xv_k\}$ if $d_{\overline{H}}(y)=1$; 
$S=\{xv_k\}$ if $d_{\overline{H}}(y) > 1$.

The proof of (A) in Claim 3 can be adapted to show the validity of (D).
Here we define
\[
H'=(H \cup
 \bigcup_{i=1}^{k-1}E(v_i\hc u_{i+1}))\setminus
 (S \cup \bigcup_{i=1}^{k-1}E(u_i \hh v_i)),
\]
where $S=\{xy,xu_k\}$ if $d_{\overline{H}}(y)=1$; 
$S=\{xu_k\}$ if $d_{\overline{H}}(y) > 1$.

In both cases, $d_{H'}(u_1)=3-1=2$ and $d_{\overline{H'}}(u_1)=2$.
It is easy to check that $H'\in \mathcal{H}(G)$ such that
$i(H')=i(H)$ and $i(\overline{H'})=i(\overline{H})-1$.  This
contradicts the choice of $H$. 
\qed

\bigskip

Now we are ready to derive contradictions from the assumption
$i(H)+i(\overline{H})>0$.

\bigskip
\noindent {\bf Case 1}  $i(H)>0$.
\bigskip

Suppose that $v_0v'\in I(H)$ with $d(v_0)\geqslant d(v')$. Let 
$\mathcal{C}(v_0)$ be the set of alternating chains of $G$ beginning 
with the vertex $v_0$.  By Claims 1 and 3, $\mathcal{C}(v_0)$ is a 
nonempty set. Let $V_{\rm I}(P)$ and $V_{\rm II}(P)$, respectively,  
be the sets of type-I vertices and type-II vertices on an alternating 
path $P \in \mathcal{C}(v_0)$. Define $V_{\rm I}=\cup\{V_{\rm I}(P) 
\mid P\in \mathcal{C}(v_0)\}$ and $V_{\rm II}=\cup\{V_{\rm II}(P) \mid 
P\in \mathcal{C}(v_0)\}$.

For any vertex $w\in V_{\rm II}$, if $x\in N_H(w)$, then either $x
\in V_{\rm I}$, or $d_H(x)=2$ and the unique vertex $y\in N_H(x)
\setminus \{w\}$ satisfies that $d_H(y)=1$ and $y\in V_{\rm I}$. Thus
\[
\sum_{z\in V_{\rm I}} d_{H}(z)\geqslant\sum_{w\in V_{\rm II}}d_H(w).
\]
Since each vertex of $V_{\rm I}$ has degree at most two in $H$, and
each vertex of $V_{\rm II}$ has degree at least two in $H$, we have
\[
2|V_{\rm I}|\geqslant \sum_{z\in V_{\rm I}} d_{H}(z)\geqslant\sum_{w
\in V_{\rm II}}d_H(w)\geqslant 2|V_{\rm II}|.
\]
Thus, $|V_{\rm I}|\geqslant |V_{\rm II}|$.

For any $z\in V_{\rm I}$, we have $d_H(z)\leqslant 2$ and $d(z)\geqslant
\Delta-1$, and hence $d_{\overline{H}}(z)\geqslant \Delta-3$. From
$d_H(v_0)=1$ and $d(v_0)=\Delta -1$, we know $d_{\overline{H}}(v_0)=\Delta
-2$. Hence,
\[
\sum_{z\in V_{\rm I}} d_{\overline{H}}(z)=d_{\overline{H}}(v_0)+
\sum_{z\in V_{\rm I} \setminus\{v_{_0}\}} d_{\overline{H}}(z)\geqslant 
|V_{\rm I}|(\Delta-3)+1.
\]
For any $w\in V_{\rm II}$, we see that $d_H(w)=3$ or $d_H(w)=
d_{\overline{H}}(w)=2$. Thus $\Delta\geqslant 6$ implies
\[
\sum_{w\in V_{\rm II}} d_{\overline{H}}(w)\leqslant |V_{\rm II}|(\Delta-3).
\]
Then $|V_{\rm I}|\geqslant |V_{\rm II}|$ implies
\[
\sum_{w\in V_{\rm II}} d_{\overline{H}}(w)< \sum_{z\in V_{\rm I}} 
d_{\overline{H}}(z).
\]
However, for $z\in V_{\rm I}$ and for each $x\in N_{\overline{H}}(z)$, 
either $x\in V_{\rm II}$, or $d_{\overline{H}}(x)=2$ and the unique 
vertex $y\in N_{\overline{H}}(x)\setminus \{w\}$ has $d_{\overline{H}}(y)=1$ 
and $y\in V_{\rm II}$. We get a contradictory consequence
\[
\sum_{w\in V_{\rm II}}d_{\overline{H}}(w) \geqslant \sum_{z\in V_{\rm I}}
d_{\overline{H}}(z).
\]

\bigskip
\noindent {\bf Case 2}  $i(\overline{H})>0$.
\bigskip

Suppose that $u_1u'\in I(\overline{H})$ with $d(u_1)\geqslant d(u')$.  
Let $\mathcal{D}(u_1)$ be the set of alternating chains of $G$
beginning with the vertex $u_1$.  By Claims 2 and 4, $\mathcal{D}(u_1)$ 
is a nonempty set. Let $V_{\rm I}(P)$ and $V_{\rm II}(P)$, respectively,  
be the sets of type-I vertices and type-II vertices on an alternating 
path $P \in \mathcal{D}(u_1)$. Define $V_{\rm I}=\cup\{V_{\rm I}(P) \mid 
P\in \mathcal{D}(u_1)\}$ and $V_{\rm II}=\cup\{V_{\rm II}(P) \mid P\in 
\mathcal{D}(u_1)\}$.

Similar to the proof of Case 1, we have that $|V_{\rm I}|\geqslant 
|V_{\rm II}|$ and
\[
|V_{\rm I}|(\Delta-3)\leqslant \sum_{z\in V_{\rm I}} d_{\overline{H}}(z)
\leqslant \sum_{w\in V_{\rm II}} d_{\overline{H}}(w).
\]
However, since $d_{\overline{H}}(u_1)=1$ and $\Delta\geqslant 6$, we get
\[
\sum_{w\in V_{\rm II}} d_{\overline{H}}(w)=d_{\overline{H}}(u_1)+\sum_{w
\in V_{\rm II}\setminus \{u_{_1}\}} d_{\overline{H}}(w)<|V_{\rm II}|(\Delta-3).
\]
A contradiction is produced. This completes the proof of Theorem \ref{p1}.
\qed

\bigskip

\noindent
{\bf Acknowledgment} This work was done while the first two authors
were visiting the third author. The support provided by the Institute
of Mathematics, Academia Sinica is greatly appreciated.

\bigskip

%%%%%%%%%%%%%%%%%%%%%%%%%%%%%%%%%%%%%%%%%%%%%%%%%%%%%%%%%%%%%%%%%%
%

\end{document}